\documentclass[12pt,a4paper]{amsart}
\usepackage{amssymb,amsmath,amsthm}

\headsep=1cm \numberwithin{equation}{section}
\theoremstyle{plain}

 \newcommand{\Z}{\mathbb{Z}}

\theoremstyle{definition}

\newtheorem{notation, definition and remark}[guess]{Notation, Definition
              and Remark}
\newtheorem{convention and remark}[guess]{Convention and Remark}
\newtheorem{notation and remark}[guess]{Notation and Remark}
\newtheorem{remark and definition}[guess]{Remark and Definition}
\newtheorem{reminder and remark}[guess]{Reminder and Remark}

\begin{document}

%%%--------------------------------------------------------------------
\title[$\Phi$-prime submodules]
 {$\phi$-prime submodules }

\author[Naser Zamani]{Naser~Zamani }

\address{faculty of science, university of mohaghegh ardabili, p.o.box 179, ardabil,
iran} \email{naserzaka@yahoo.com, zamanin@uma.ac.ir}

%\thanks{This work was completed with the support of an Izaak
% Walton Killam Memorial Scholarship.}

%\thanks{The author was also supported in part by the Research
 %Council of Slovenia.}



%\date{February 15, 1995 and, in revised form, July 6, 1995.}

%\dedicatory{}

%\commby{Daniel J. Rudolph}

%%% ----------------------------------------------------------------------

\begin{abstract}
Let $R$ be a commutative ring with non-zero identity and $M$ be a
unitary $R$-module. Let $\mathcal{S}(M)$ be the set of all
submodules of $M$, and $\phi:\mathcal{S}(M)\rightarrow
\mathcal{S}(M)\cup \{\emptyset\}$ be a function. We say that a
proper submodule $P$ of $M$ is a prime submodule relative to $\phi$
or $\phi$-prime submodule if $a\in R$, $x\in M$ with $ax\in
P\setminus \phi(P)$ implies that $a\in(P:_RM)$ or $x\in P$. So if we
take $\phi(N)=\emptyset$ for each $N\in\mathcal{S}(M)$, then a
$\phi$-prime submodule is exactly a prime submodule. Also if we
consider $\phi(N)=\{0\}$  for each submodule $N$ of $M$, then in
this case a $\phi$-prime submodule will be called a weak prime
submodule.
 Some of the properties of this concept will be
investigated. Some characterizations of $\phi$-prime submodules will
be given, and we show that under some assumptions prime submodules
and $\phi_1$-prime submodules coincide.
\end{abstract}

%%% ----------------------------------------------------------------------
\maketitle
\section{introduction} Throughout $R$ is a commutative ring with
nonzero identity and $M$ is a unitary $R$-module. Prime ideals play
an essential role in ring theory. One of the natural generalizations
of prime ideals which have attracted the interest of several authors
in the last two decades is the notion of prime submodule, (see for
example [2-6]). These have led to more information on the structure
of the $R$-module $M$. For an ideal $I$ of $R$ and a submodule $N$
of $M$ let $\sqrt{I}$ denote the radical of $I$, and $(N:_RM)=\{r\in
R: rM\subseteq N\}$ which is clearly a submodule of $M$. We say that
$N$ is a radical submodule of $M$ if $\sqrt{(N:_RM)}=(N:_RM)$. Then
a proper submodule $P$ of $M$ is called a prime submodule if $r\in
R$ and $x\in M$ with $rx\in P$ implies that $r\in (P:_RM)$ or $x\in
P$. It is easy to see that $P$ is a prime submodule of $M$ if and
only if $(P:_RM)$ is a prime ideal of $R$ and the $R/(P:_RM)$-module
$M/P$ is torsion free (the $R$-module $X$ is said to be torsion free
if the annihilator
 of any nonzero element of $X$ is zero). By
restricting where $rx$ lies we can generalize this definition. A
submodule $P\neq M$ is said to be \textit{weak prime submodule} of
$M$ if $r\in R$ and $x\in M$, $0\neq rx\in P$ gives that $r\in
(P:_RM)$ or $x\in P$. We will say that $P\neq M$ is almost prime
submodule if $r\in R$ and $x\in M$ with $rx\in P\setminus(P:_RM)P$
implies that $r\in (P:_RM)$ or $x\in P$. So any prime submodule is
weak prime and any weak prime submodule is an almost prime
submodule. Let $\mathcal{S}(M)$ be the set of all submodules of $M$
and $\phi: \mathcal{S}(M)\rightarrow
\mathcal{S}(M)\cup\{\emptyset\}$ be a function. A proper submodule
$P$ of $M$ is said to be a $\phi$-prime submodule if $r\in R$ and
$x\in M$, $rx\in P\setminus\phi(P)$ implies that $r\in (P:_RM)$ or
$x\in P$. Since $P\setminus \phi(P)=P\setminus(P\cap\phi(P))$, so
without loss of generality, throughout this paper we will consider
$\phi(P)\subseteq P$. In the rest of the paper we use the following
functions $\phi:\mathcal{S}(M)\rightarrow
 \mathcal{S}(M)\cup\{\emptyset\}$.
\[\phi_\emptyset(N)=\emptyset, \quad \quad \quad \quad \quad \quad \forall \quad N\in \mathcal{S}(M),\]
 \[\phi_0(N)=\{0\}, \quad \quad \quad \quad \quad \quad \forall \quad N \in \mathcal{S}(M),\]
 \[\phi_1(N)=(N:_RM)N,  \quad \quad  \quad\forall \quad N \in \mathcal{S}(M),\]
 \[\phi_2(N)=(N:_RM)^2N,\quad  \quad\quad \forall \quad N \in \mathcal{S}(M),\]
 \[\phi_\omega(N)=\cap_{i=1}^\infty (N:_RM)^iN, \quad \quad \forall \quad N \in \mathcal{S}(M).\]

 Then it is clear that $\phi_{\emptyset}$, $\phi_0$-prime submodules are prime,
 weak prime submodules respectively. Evidently for any submodule and every
 positive integer $n$, we have the following implications:
 \[\textrm{prime}\Rightarrow\phi_\omega-\textrm{prime}\Rightarrow\phi_n-\textrm{prime}
 \Rightarrow\phi_{n-1}-\textrm{prime}.\]For functions
   $\phi, \psi: \mathcal{S}(M)\rightarrow
  \mathcal{S}(M)\cup\{\emptyset\}$, we write $\phi\leq \psi$ if $\phi(N)\subseteq \psi(N)$
   for each $N\in\mathcal{S}(M)$. So whenever $\phi\leq \psi$,
    any $\phi-$prime submodule is $\psi$-prime.\\
    In this paper, among other
  results concerning the properties of $\phi$-prime submodules, some
  characterizations of this notion will be investigated. Some of the
  results in this paper inspired from \cite{And-Bat}.
 % -----------------------------------------------------------------------------
\section{results}
The following Theorem asserts that under some conditions
$\phi$-prime submodules are prime.\\

 \textbf{Theorem 2.1.} \textit{Let $R$ be a commutative ring and $M$ be an
 $R$-module. Let $\phi: \mathcal{S}(M)\rightarrow
 \mathcal{S}(M)\cup\{\emptyset\}$ be a function and $P$ be a
 $\phi$-prime submodule of $M$ such that $(P:_RM)P\nsubseteq \phi(P)$.
 Then $P$ is a prime submodule of $M$.}
\begin{proof} Let $a\in R$ and
$x\in M$ be such that $ax\in P$. If $ax\notin \phi(P)$, then since $P$
is $\phi$-prime, we have $a\in (P:_RM)$ or $x\in P$.\\
So let $ax\in \phi(P)$. In this case we may assume that $aP\subseteq
\phi(P)$. For, let $aP\nsubseteq \phi(P)$. Then there exists $p\in
P$ such tat $ap\notin \phi(P)$, so that $a(x+p)\in
P\setminus\phi(P)$. Therefore $a\in (P:_RM)$ or $x+p\in P$ and hence
$a\in (P:_RM)$ or $x\in P$. Second we may assume that
$(P:_RM)x\subseteq \phi(P)$. If this is not the case, there exists
$u\in (P:_RM)$ such that $ux\notin \phi(P)$ and so $(a+u)x\in
P\setminus\phi(P)$. Since $P$ is a $\phi$ -prime submodule, we have
$a+u\in (P:_RM)$ or $x\in P$. So $a\in (P:_RM)$ or $x\in P$. Now
since $(P:_RM)P\nsubseteq \phi(P)$, there exist $r\in (P:_RM)$ and
$p\in P$ such that $rp\notin\phi(P)$. So $(a+r)(x+p)\in
P\setminus\phi(P)$, and hence $a+r\in(P:_RM)$ or $x+p\in P$.
Therefore $a\in (P:_RM)$ or $x\in P$ and the proof is complete.
\end{proof}

\textbf{Corollary 2.2.}\textit{ Let $P$ be a weak prime submodule of
$M$ such that $(P:_RM)P\neq 0$. Then $P$ is a prime submodule of
$M$.}
\begin{proof} In the above Theorem set $\phi=\phi_0$.
\end{proof}
 \textbf{Corollary 2.3.} \textit{Let $P$ be a
$\phi$-prime submodule of $M$ such that $\phi(P)\subseteq
(P:_RM)^2P$. Then for each $a\in R$ and $x\in M$, $ax\in
P\setminus\cap_{i=1}^{\infty}(P:_RM)^iP$ implies that $a\in(P:_RM)$
or $x\in P$. In the other word $P$ is $\phi_\omega$-prime.}
\begin{proof} If $P$ is a prime submodule of $M$, then the result is clear.
 So suppose that $P$ is not a prime submodule of $M$.
 Then by Theorem 2.1 we have $(P:_RM)P\subseteq \phi(P)\subseteq (P:_RM)^2P\subseteq(P:_RM)P$,
  that is $\phi(P)=(P:_RM)P=(P:_RM)^2P$.
 Hence $\phi(P)=(P:_RM)^iP$ for all $i\geq 1$ and the result follows.
\end{proof}
\textbf{Corollary 2.4.} \textit{Let $M$ be an $R$ module and $P$ be
a $\phi$-prime submodule of $M$. Then} $(P:_RM)\subseteq
\sqrt{(\phi(P):_RM)}$ \textit{or}
$\sqrt{(\phi(P):_RM)}\subseteq(P:_RM)$. \textit{If}
$(P:_RM)\varsubsetneq \sqrt{(\phi(P):_RM)}$ \textit{then $P$ is not
a prime submodule of $M$; while if}
$\sqrt{(\phi(P):_RM)}\varsubsetneq (P:_RM)$, \textit{then $P$ is a
prime submodule of $M$. If $\phi(P)$ is a radical submodule of $M$,
either $(P:_RM)=(\phi(P):_RM)$ or $P$ is a prime submodule of $M$.}
\begin{proof} If $P$ is not a prime submodule of $M$, then by Theorem 2.1,
we have $(P:_RM)P\subseteq \phi(P)$. Hence
$\sqrt{(P:_RM)^2}\subseteq \sqrt{((P:_RM)P:_RM)}\subseteq
 \sqrt{(\phi(P):_RM)}$. So $(P:_RM)\subseteq \sqrt{(\phi(P):_RM)}$. If $P$ is a prime
  submodule of $M$, then $\sqrt{(\phi(P):_RM)}\subseteq \sqrt{(P:_RM)}=(P:_RM)$
  (note that we may assume that $\phi(P)\subseteq P$), and all the claims of the corollary follows.
\end{proof}
\textbf{Remark A.} Suppose that $P$ is a $\phi$-prime submodule of
$M$ such
 that $\phi(P)\subseteq (P:_RM)P$ (resp. $\phi(P)\subseteq (P:_RM)^2P$) and that
   $P$ is not a prime submodule. Then by Theorem 2.1, we have $\phi(P)=(P:_RM)P$
   (resp. $\phi(P)=(P:_RM)^2P$). In particular if $P$ is a weak prime
    (resp. $\phi_2$-prime) submodule but not prime submodule then $(P:_RM)P=0$
     (resp. $(P:_RM)P=(P:_RM)^2P$).\\

Let $R_1, R_2$ be two commutative rings with identity. Let $M_1$ and
$M_2$ be $R_1$ and $R_2$-module respectively and put $R=R_1\times
R_2$. Then $M=M_1\times M_2$ is an $R$-module and each submodule of
$M$ is of the form $N=N_1\times N_2$ for some submodule $N_1$ of
$M_1$ and $N_2$ of $M_2$. Furthermore $N=N_1\times N_2$ is a prime
submodule of $M$ if and only if $N=P_1\times M_2$ or $N=M_1\times
P_2$ for some prime submodule $P_1$ of $M_1$ and $P_2$ of $M_2$. (
In fact, to see the nontrivial direction, let $N_1\times N_2$ be a
prime submodule of $M_1\times M_2$. Then either $N_1$ must be a
prime submodule of $M_1$ or $N_2$ must be a prime submodule of
$M_2$. Now $(N_1:_{R_1}M_1)\times (N_2:_{R_2}M_2)=(N_1\times
N_2:_RM_1\times M_2)$ is a prime ideal of $R=R_1\times R_2$. So
either $(N_1:_{R_1}M_1)=R_1$ or $(N_2:_{R_2}M_2)=R_2$, which means
that either $N_1=M_1$ or $N_2=M_2$ and the claim follows.) If $P_1$
is a weak
 prime submodule of $M_1$, then $P_1\times M_2$ need not be a weak prime submodule of $M$.
  Indeed $P_1\times M_2$ is a weak prime submodule of $M$ if and only if $P_1\times M_2$
  is a prime
submodule of $M_1\times M_2$. To see the nontrivial direction, let
$P_1\times M_2$ be a weak prime submodule of $M_1\times M_2$. Let
$r_1\in R_1$, $x_1\in M_1$ with $r_1x_1\in P_1$. Let $0\neq x_2\in
M_2$. Then $(r_1, 1)(x_1,x_2)=(r_1x_1,x_2)\in P_1\times
M_2\setminus\{(0,0)\}$. By assumption this gives that
$(r_1,1)\in(P_1\times M_2:_{R_1\times R_2}M_1\times
M_2)=(P_1:_{R_1}M_1)\times R_2$ or $(x_1, x_2)\in P_1\times M_2$,
that is $r_1\in(P_1:_{R_1}M_1)$ or $x_1\in P_1$. Therefore $P_1$ is
a prime submodule of $M_1$ and
hence $P_1\times M_2$ is a prime submodule of $M_1\times M_2$.\\
However if $P_1$ is a weak prime submodule of $M_1$, then $P_1\times
M_2$ is a $\phi$-prime submodule if $\{0\}\times M_2\subseteq
\phi(P_1\times M_2)$.

    To see this, we have
    $P_1\times M_2\setminus\phi(P_1\times M_2)
    \subseteq P_1\times M_2\setminus\{0\}\times M_2
    =(P_1\setminus\{0\})\times M_2$. Now let $(r_1,r_2)(x_1, x_2)
    =(r_1x_1,r_2x_2)\in P_1\times M_2\setminus\phi(P_1\times M_2)$.
    Then $r_1x_1\in P_1\setminus\{0\}$ and by the assumption on $P_1$ we have
     $r_1\in (P_1:_{R_1}M_1)$ or $x_1\in P_1$. This gives that $(r_1,r_2)\in (P_1:_TM_1)
     \times R_2 =(P_1\times M_2:_{R_1\times R_2}M_1\times M_2)$ or
     $(x_1, x_2)\in P_1\times M_2$.
     Therefore $P_1\times M_2$ is a $\phi$-prime submodule of $M_1\times M_2$.\\

\textbf{Corollary 2.5.} \textit{Let $R_1$ and $R_2$ be two
commutative rings,
 $M_1$ and $M_2$ be $R_1$ and $R_2$-modules respectively. Let $M=M_1\times M_2$
 and $\phi:\mathcal{S}(M)\rightarrow \mathcal{S}(M)\cup\{\emptyset\}$ be a function.
 Suppose that $P_1$ is a
 weak prime submodule of $M_1$ such that $\{0\}\times M_2\subseteq \phi(P_1\times
 M_2)$.
 %(e.g., $\phi_\omega \leq \phi$).
 Then $P_1\times M_2$ is a $\phi$-prime submodule of $M_1\times M_2$.}\\

 \textbf{Proposition 2.6.} \textit{With the same notations as in Corollary 2.5,
 Let $\phi$ be a function such that $\phi_\omega\leq
 \phi$. Then for any weak prime submodule  $P_1$ of $M_1$, $P_1\times M_2$ is a
 $\phi$-prime submodule of $M_1\times M_2$.}
\begin{proof} If $P_1$ is a prime submodule of $M_1$,
 then $P_1\times M_2$ is prime and so a $\phi$-prime submodule of $M_1\times M_2$.
 Suppose that $P_1$ is not a prime submodule of $M_1$. Then by Remark A, we
 have $(P_1:_{R_1}M_1)P_1=0$. This gives
   that
   \[(P_1\times M_2:_{R_1\times R_2}M_1\times M_2)^i(P_1\times M_2)=[(P_1:_{R_1}M_1)^iP_1]\times M_2=0\times
   M_2,\]
    for all $i\geq 1$ and hence we have
    \[0\times M_2=\cap_{i=1}^\infty(P_1\times M_2:_{R_1\times R_2}M_1\times M_2)^i
    (P_1\times M_2)=\phi_\omega(P_1\times M_2)
      \subseteq \phi(P_1\times M_2),\] and the result follows by the
   above Corollary.
   \end{proof}

\textbf{Theorem 2.7.} \textit{Let $M$ be an $R$-module and $0\neq
x\in M$  such that $Rx\neq M$ and $(0:_Rx)=0$. If $Rx$ is not a
prime submodule of $M$, then $Rx$ is not a $\phi_1$-prime submodule
of $M$. }
\begin{proof} Since $Rx$ is not a prime submodule of $M$, there exists $a\in R$,
 $y\in M$ such that, $a\notin (Rx:_RM)$, $y\notin Rx$, but $ay\in Rx$.
 If $ay\notin(Rx:_RM)x$, then by our definition $Rx$ is not a $\phi_1$-prime submodule.
 So let $ay\in (Rx:_RM)x$. We have $y+x\notin Rx$ and $a(y+x)\in Rx$.
  If $a(y+x)\notin (Rx:_RM)x$, then again by our definition $Rx$ is not
   $\phi_1$-prime submodule.
  So let  $a(y+x)\in (Rx:_RM)x$, then $ax\in(Rx:_RM)x$, which gives that $ax=rx$ for
   some $r\in (Rx:_RM)$. Since $(0:_Rx)=0$ this gives
    that $a=r\in (Rx:_RM)$ which  contradicts with our assumption.
\end{proof}
\textbf{Corollary 2.8.} \textit{Let $x$ be a nonzero element of an
$R$-module $M$ such that $(0:_Rx)=0$ and that $Rx\neq M$. Then $Rx$
is a prime
 submodule of $M$ if and only if $Rx$ is a $\phi_1$-prime submodule of $M$.}\\

\textbf{Proposition 2.9.} \textit{ Let $P$ be a $\phi_1$-prime
submodule of $M$. Then the following holds:}\\
(i) \textit{ If $a$ is a zero divisor in $M/P$, then $aP\subseteq(P:_RM)P$.}\\
(ii) \textit{ Let $J$ be an ideal of $R$ such that $(P:_RM)\subseteq J$ and $J\subseteq Z_R(M/P)$, then $JP=(P:_RM)P$.}
\begin{proof}
(i) By assumption, there exists $x\in M\setminus P$ such that $ax\in
P$. If $a\in (P:_RM)$ then clearly $aP\subseteq (P:_RM)P$. So let
$a\notin (P:_RM)$.
 Since $P$ is a $\phi_1$-prime submodule of $M$, we must
 have $ax\in (P:_RM)P$.  Now for any $y\in P$; $y+x\notin P$ and $a(y+x)\in P$.
  Hence as $P$ is $\phi_1$-prime $a(y+x)\in (P:_RM)P$,
   which gives that $ay\in (P:_RM)P$. So $aP\subseteq (P:_RM)P$ and the result follows.

(ii) This follows from (i).
\end{proof}

\textbf{Theorem 2.10.} \textit{Let $M$ be an $R$-module and let $a$
be an element of $R$ such that $aM\neq M$. Suppose $(0:_Ma)\subseteq
aM$. Then $aM$ is a $\phi_1$-prime submodule of $M$ if and only if
it is a prime submodule of $M$.}
\begin{proof}
The direction $\Leftarrow$ is clear. So we prove $\Rightarrow$.
Let $b\in R$, $x\in M$ such that $bx\in aM$. We show that $b\in (aM:_RM)$
or $x\in aM$. If $bx\notin (aM:_RM)aM$, then $b\in (aM:_RM)$ or $x\in aM$,
 since $aM$ is a $\phi_1$-prime submodule. So suppose $bx\in (aM:_RM)aM$.
 Now $(b+a)x\in aM$. If $(b+a)x\notin (aM:_RM)aM$, then, since $aM$ is a $\phi_1$
 prime submodule, $b+a\in (aM:_RM)$ or $x\in aM$ which give that $b\in (aM:_RM)$
  or $x\in aM$. So assume that $(b+a)x\in (aM:_RM)aM$. Then $bx\in(aM:_RM)aM$ gives that $ax\in(aM:_RM)aM$. Hence there exists $y\in (aM:_RM)M$ such that $ax=ay$ and so $x-y\in (0:_Ma)$. This gives that $x\in (aM:_RM)M+(0:_Ma)\subseteq aM+(0:_Ma)\subseteq aM$, and the result follows.
\end{proof}
In the next theorem we give several characterizations of $\phi$-prime submodules.\\

\textbf{Theorem 2.11.} \textit{ Let $P$ be a proper submodule of $M$ and let
 $\phi:\mathcal{S}(M)\rightarrow \mathcal{S}(M)\cup \{\emptyset\}$ be a function.
  Then the following are equivalen:}\\
(i)\textit{ $P$ is a $\phi$-prime submodule of $M$;}\\
(ii)\textit{ for $x\in M\setminus P$, $(P:_Rx)=(P:_RM)\cup (\phi(P):_Rx)$;}\\
(iii)\textit{ for $x\in M\setminus P$, $(P:_Rx)=(P:_RM)$ or $(P:_Rx)=(\phi(P):_Rx)$;}\\
(iv)\textit{ for any ideal $I$ of $R$ and any submodule $L$ of $M$, if $IL\subseteq P$ and $IL\nsubseteq \phi(P)$, then $I\subseteq (P:_RM)$ or $L\subseteq P$.}
\begin{proof}(i)$\Rightarrow $(ii). Let $x\in M\setminus P$ and
$a\in (P:_Rx)\setminus(\phi(P):_Rx)$. Then $ax\in
P\setminus\phi(P)$. Since $P$ is a $\phi$-prime submodule of $M$, So
$a\in (P:_RM)$. As we may assume that
$\phi(P)\subseteq P$, the other inclusion always holds.\\
(ii)$\Rightarrow $ (iii). If a subgroup is the union of two
subgroups, it is equal to one of them.\\
(iii)$\Rightarrow $(iv). Let $I$ be an ideal of $R$ and $L$ be a submodule of $M$
such that $IL\subseteq P$. Suppose $I\nsubseteq (P:_RM)$ and $L\nsubseteq P$.
We show that $IL\subseteq \phi(P)$. Let $a\in I$ and $x\in L$.
 First let $a\notin (P:_RM)$. Then, since $ax\in P$, we have $(P:_Rx)\neq (P:_RM)$.
 Hence by our assumption $(P:_Rx)=(\phi(P):_Rx)$. So $ax\in \phi(P)$.
  Now assume that $a\in I\cap (P:_RM)$. Let $u\in I\setminus(P:_RM)$.
  Then $a+u\in I\setminus(P:_RM)$. So by the first case,
  for each $x\in L$ we have $ux \in \phi(P)$ and $(a+u)x\in \phi(P)$.
  This gives that $ax\in\phi(P)$. Thus in any case $ax \in \phi(P)$.
   Therefore $IL\subseteq \phi(P)$. \\
(iv)$\Rightarrow$ (i). Let $ax\in P\setminus\phi(P)$. By considering
the ideal $(a)$ and the submodule $(x)$, the result follows.
\end{proof}

Let $S$ be a multiplicatively close subset of $R$.
 Then by \cite[9.11 (v)]{Sharp} each submodule of $S^{-1}M$ is of
 the form $S^{-1}N$ for some submodule $N$ of $M$.
 Also it is well known that there is a one to one correspondence between
 the set of all prime submodules $P$ of $M$ with $(P:_RM)\cap S=\emptyset$
  and the set of all prime submodules of $S^{-1}M$, given by $P\rightarrow
  S^{-1}P$ (see \cite[Theorem 3.4]{Mo-Sally}). Furthermore  it is easy to see that if $P$
   is a weak prime submodule of $M$ with $S^{-1}P\neq S^{-1}M$, then
   $S^{-1}P$ is a weak prime submodule of $S^{-1}M$. This fact remains
   true for $\phi_1$-prime submodules $P$ of $M$ with $S^{-1}P\neq S^{-1}M$. In the next theorem we want to generalize this fact
   for $\phi$-prime submodules. In the following for a submodule
    $N$ of $M$ we put $N(S)=\{x\in M:
  \exists s\in S, sx\in N\}.$ Then $N(S)$ is a submodule of $M$
  containing $N$ and $S^{-1}(N(S))=S^{-1}N$. Let
  $\phi:\mathcal{S}(M)\rightarrow \mathcal{S}(M)\cup\{\emptyset\}$
  be a function. We define
  $(S^{-1}\phi):\mathcal{S}(S^{-1}M)\rightarrow
  \mathcal{S}(S^{-1}M)\cup\{\emptyset\}$ by $(S^{-1}\phi)(S^{-1}N)=
  S^{-1}(\phi(N(S)))$ if $\phi(N(S))\neq \emptyset$ and
  $(S^{-1}\phi)(S^{-1}N)=\emptyset$ if $\phi(N(S))= \emptyset$.
  Since dealing with prime submodules $P$ we can always assume
  that $\phi(P)\subseteq P$, so there is no loss of generality in
  assuming that $\phi(N)\subseteq N$, and hence
  $(S^{-1}\phi)(S^{-1}N)\subseteq S^{-1}N$. Also we note that
  $(S^{-1}\phi_\emptyset)=\phi_\emptyset$, $(S^{-1}\phi_0)=\phi_0$,
  and whenever $M$ is finitely generated
  $(S^{-1}\phi_i)=\phi_i$ for $i=1,2$. In the next theorem we show
  that if $S^{-1}(\phi(N))\subseteq (S^{-1}\phi)(S^{-1}N)$, then
  $\phi$-primeness of $P$ together with $S^{-1}P\neq S^{-1}M$ imply  that
  $S^{-1}P$ is $(S^{-1}\phi)$-prime.\\
  For a submodule $L$ of $M$, we define
  $\phi_L:\mathcal{S}(M/L)\rightarrow
  \mathcal{S}(M/L)\cup\{\emptyset\}$ by $\phi_L(N/L)=(\phi(N)+L)/L$
  for $N\supseteq L$ and $\emptyset$ for $\phi(N)=\emptyset$.\\

  \textbf{Theorem 2.12.}\textit{ Let $M$ be an $R$-module and let
  $\phi:\mathcal{S}(M)\rightarrow \mathcal{S}(M)\cup \{\emptyset\}$.
  Let $P$ be a $\phi$-prime submodule of $M$.}\\
  (i)\textit{ If $L\subseteq P$ is a submodule of $M$, then $P/L$ is
  a  $\phi_L$-prime submodule of $M/L$.}\\
  (ii)\textit{ Suppose that $S$ is a multiplicatively closed subset
  of $R$ such that $S^{-1}P\neq S^{-1}M$ and
  $S^{-1}(\phi(P))\subseteq (S^{-1}\phi)(S^{-1}P)$. Then $S^{-1}P$
  is an $(S^{-1}\phi)$-prime submodule of $S^{-1}M$. Furthermore  if $S^{-1}P\neq
  S^{-1}((\phi(P))$, then P(S)=P.}
  \begin{proof}
  (i) Let $a\in R$ and $\bar{x}\in M/L$ with
  $a\bar{x}\in P/L\setminus\phi_L(P/L)$, where $\bar{x}=x+L$, for some $x\in M$.
   By the definition of $\phi_L$, this gives that
  $ax\in P\setminus(\phi(P)+L)$. So we have $ax\in P\setminus\phi(P)$,
  which gives that $a\in (P:_RM)$ or
   $x\in P$. Therefor $a\in (P/L:_RM/L)$ or $x\in P$ and so $P/L$ is $\phi_L$-prime
   submodule.\\
   (ii) Let $a/ s \in S^{-1}R$ and $x/t \in S^{-1}M$ with
    $ax/st\in S^{-1}P\setminus(S^{-1}\phi)(S^{-1}P)$. Then by our assumption
     $ax/st\in S^{-1}P\setminus S^{-1}(\phi(P))$.
    Therefore there exists
    $u\in S$ such that $uax\in P\setminus\phi(P)$, (note that for each $v\in
    S$, $vax\notin \phi(P))$. Since $P$ is $\phi$-prime and
    $(P:_RM)\cap S=\emptyset$, this gives that $ax\in P\setminus\phi(P)$ and
    so $a\in (P:_RM)$ or $x\in P$. Therefore $a/s\in
    S^{-1}(P:_RM)\subseteq (S^{-1}P:_{S^{-1}R}S^{-1}M)$ or $x/t\in
    S^{-1}P$. Hence $S^{-1}P$ is an $(S^{-1}\phi)$-prime submodule of
    $S^{-1}M$.\\
    To prove the last part of the theorem, let $x\in P(S)$. Then
    there exists $s\in S$ such that $sx\in P$. If $sx\notin \phi(P)$,
    then $x\in P$. If $sx\in \phi(P)$, then $x\in\phi(P)(S)$. So
    $P(S)=P\cup (\phi(P)(S))$. Hence $P(S)=P$ or $P(S)=(\phi(P)(S))$. If
    the second holds, then we must have
    $S^{-1}P=S^{-1}P(S)=S^{-1}(\phi(P)(S))=S^{-1}(\phi(P))$, which is
    not the case.  So $P(S)=P$ and the proof is complete.
  \end{proof}

Let $S^{-1}P$ be an $(S^{-1}\phi)$-prime submodule of $S^{-1}M$.
Then evidently $(P:_RM)\cap S=\emptyset$. In general we don't know
under
 what conditions $P$ is a $\phi$-prime submodule of $M$. Even in the
 case $\phi=\phi_0, \phi_1$ and $\phi_2$ we could not answer this
 question.\\
 As we mentioned previously, for two commutative rings $R_1$ an $R_2$ and
 two modules $M_1$ and $M_2$ over $R_1$ and $R_2$ respectively, the
 prime submodules of the $R=R_1\times R_2$ module $M=M_1\times M_2$
 are in the form $P_1\times
 M_2$ or $M_1\times P_2$ where $P_1$ is a prime submodule of $M_1$
 and $P_2$ is a prime  submodule of $M_2$. This is not true for
 correspondence $\phi$-prime submodules in general. For example if $P_1$ is a $\phi_0$-prime
 submodule of $M_1$, then $P_1\times M_2$ is not necessarily a
 $\phi_0$-prime submodule of $M_1\times M_2$. To be more specific
 let $R_1=R_2=M_1=M_2=\Z_6$, and suppose $P_1=\{0\}$. Then evidently
 $P_1$ is a $\phi_0$-prime submodule of $M_1$. However,
 $(2,1)(3,1)\in P_1\times M_2$, and $(3,1)\notin P_1\times M_2$. Also as
 $(2,1)(2,1)\notin P_1\times M_2$, $(2,1)M\nsubseteq P_1\times M_2$.
 However in this direction we have the
  following result. \\

  \textbf{Theorem 2.13.} \textit{Let the notation be as in the above
  paragraph. Let $\psi_i:\mathcal{S}(M_i)\rightarrow
  \mathcal{S}(M_i)\cup\{\emptyset\}$. Let $\phi=\psi_1\times
  \psi_2$. Then each of the following types
  are $\phi$-prime submodules of $M_1\times M_2$:\\}
 (i) \textit{$N_1\times N_2$ where $N_i$ is a proper submodule of
 $M_i$, with $\psi_i(N_i)=N_i$.}\\
(ii) \textit{$P_1\times M_2$ where $P_1$ is a prime
 submodule of $M_1$.}\\
 (iii) \textit{$P_1\times M_2$ where $P_1$ is a $\psi_1$-prime
 submodule of $M_1$ and $\psi_2(M_2)= M_2$.}\\
(iv) \textit{$M_1\times P_2$ where $P_2$ is a prime submodule of
$M_2$.}\\
 (v) \textit{$M_1\times P_2$ where $P_2$ is a $\psi_2$-prime
submodule of $M_2$  and $\psi_1(M_1)= M_1$.}
 \begin{proof}
 (i) is clear, since $N_1\times N_2\setminus\phi(N_1\times
 N_2)=\emptyset$. \\
  (ii) If $P_1$ is a prime submodule of $M_1$, then $P_1\times
 M_2$ as a prime submodule of $M_1\times M_2$ is $\phi$-prime.\\
 (iii)
  Let $P_1$ be a  $\psi_1$-prime submodule of $M_1$ and
 $\psi_2(M_2)=M_2$. Let $(r_1, r_2)\in R$ and $(x_1, x_2)\in M$ be
 such that $(r_1, r_2)(x_1, x_2)=(r_1x_1, r_2x_2)\in P_1\times M_2\setminus\phi(P_1\times
 M_2)=P_1\times M_2\setminus\psi_1(P_1)\times\psi_2(M_2)=P_1\times
 M_2\setminus\psi_1(P_1)\times M_2=(P_1\setminus\psi_1(P_1))\times M_2.$ So
$r_1\in(P_1:_{R_1}M_1)$ or $x_1\in P_1$. Therefore  $(r_1,
 r_2)\in(P_1:_{R_1}M_1)\times
 R_2=(P_1\times M_2:_{R_1\times R_2}M_1\times M_2)$ or $(x_1,x_2)\in
 P_1\times M_2)$. So $P_1\times M_2$ is a $\phi$-prime submodule of
$M_1\times M_2$.\\
  Parts  (iv), (v) are proved similar to (ii), (iii) respectively.
 \end{proof}
 A question arises  here  is that if any prime submodule of $M$ has one of the above forms.
 As it has been shown in \cite[Theorem 16]{And-Bat} this is true for the ideal and the ring case.
 But we were not able to prove the similar result in the module
 case.\\
\textbf{Acknowledgement.} I would like to thank the referee for
his/her valuable comments and suggestions.

  \bibliographystyle{plain}
      
       \end{document}